\documentclass[graybox, envcountsame, envcountsect]{svmult-fgs}
\usepackage{newtxtext,newtxmath}
\usepackage{helvet}
\usepackage{courier}
\usepackage{graphicx}
\usepackage[bottom]{footmisc}
\makeindex

\begin{document}

\title*{The random conductance model with heavy tails on nested fractal graphs}
\author{David A.\ Croydon}
\institute{David A.\ Croydon \at Research Institute for Mathematical Sciences, Kyoto University, Kyoto 606-8502, Japan, \email{croydon@kurims.kyoto-u.ac.jp}}
\maketitle

\abstract{Recently, Kigami's resistance form framework has been applied to provide a general approach for deriving the scaling limits of random walks on graphs with a fractal scaling limit \cite{C,CHK}. As an illustrative example, this article describes an application to the random conductance model with heavy tails on nested fractal graphs.
\\
\subclassname{Primary: 28A80; Secondary: 60K37}}

\section{Introduction}

One of the early motivations for the study of stochastic processes on fractals came from physics, where there was an interest in understanding the dynamical properties of disordered media. Specifically, certain examples of the latter were modelled by critical percolation, which is believed to exhibit large scale fractal structure. (See \cite{AH} for background.) The initial response from the mathematics community was to construct Brownian motion on idealised fractals, such as the Sierpi\'nski gasket \cite{Gold,Kusu}. Since then, the technology has developed to the point where it can engage with some of the original questions about critical percolation. For instance, recent work in this direction underlines that the notion of a resistance form, as introduced by Kigami to provide a broad framework for studying analysis on fractals \cite{kig1,Kig}, is useful for understanding the scaling limits of various models of random walks on random graphs in critical regimes \cite{C,CHK}. We highlight that resistance forms are only really applicable in low-dimensional settings, with the stochastic processes constructed from them typically being point recurrent (note that in the case of the standard Brownian motion on $\mathbb{R}^d$, the latter property holds only when $d=1$, and this is indeed the only dimension in which the Brownian motion can be described by a resistance form). A brief introductory survey of the work of  \cite{C,CHK} already appears in \cite{C2}, where a number of applications to random graphs are listed (see also \cite{Andrio,Archer} for some further ones that have appeared more recently), and a conjecture for critical percolation is made. Here, the aim will be to introduce the general resistance form results of \cite{C,CHK} specifically to an audience that has some familiarity with analysis on self-similar fractals by presenting in detail an example from \cite{CHK} which is of interest in its own right: the random conductance model with heavy tails on nested fractal graphs.

\begin{figure}[t]
\begin{center}
\includegraphics[width=\textwidth]{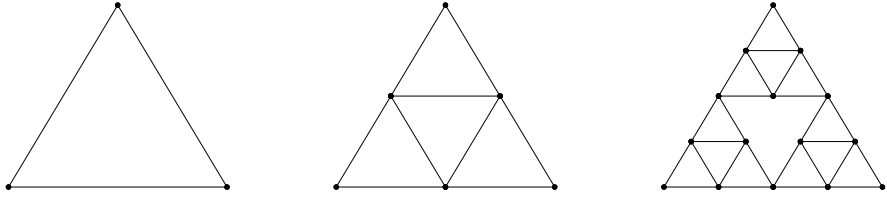}
\end{center}
\caption{The Sierpi\'nski gasket graphs $G_0$, $G_1$, $G_2$.}\label{sg}
\end{figure}

The nested fractals were originally introduced in \cite{lind}, and are a class of self-similar fractals that are finitely-ramified, embedded into Euclidean space and admit a high degree of symmetry. In the next section we will introduce sequences of graphs associated with nested fractals, but to keep the presentation concise here, we focus for the moment on a concrete example of a nested fractal, the Sierpi\'nski gasket in two dimensions. Let $V_0:=\{x_0,x_1,x_2\}\subseteq\mathbb{R}^2$ consist of the vertices of an equilateral triangle of side length 1. Write $\psi_i(x):=|x+x_i|/2$ for $i=0,1,2$. Then there exists a unique compact set $F$ such that $F=\cup_{i=0}^2\psi_i(F)$; this is the Sierpi\'nski gasket. We define the associated Sierpi\'nski gasket graphs $(G_n)_{n\geq 0}$ by setting the vertex set $V(G_n):=V_n$, where $V_n:=\cup_{i=0}^2\psi_i(V_{n-1})$ for $n\geq 1$, (note that $V_0$ was already defined,) and defining the edge set $E(G_n)$ to be the collection of pairs of elements of $V_n$ at a Euclidean distance $2^{-n}$ apart. (The first three graphs in this sequence are shown in Figure \ref{sg}.) For each $n$, we associate a stochastic process $X^n=(X^n_t)_{t\geq 0}$ by supposing $X^n$ is the continuous time Markov chain that has exponential holding times of unit mean, and at jump times moves to a neighbour of the current location with uniform probability amongst the possibilities. If we moreover assume that $X^n_0=x_0$ for each $n$, then, from the seminal early works in the area \cite{BP,Gold,Kusu,lind} it is known that
\begin{equation}\label{rwscaling}
\left(X^n_{5^nt}\right)_{t\geq 0}\rightarrow \left(X^{SG}_t\right)_{t\geq 0}
\end{equation}
in distribution in $D([0,\infty),\mathbb{R}^2)$ (that is, the space of cadlag processes on $\mathbb{R}^2$, i.e.\ those that are right-continuous and have left-hand limits, equipped with the usual Skorohod $J_1$-topology -- for elementary introductions to this framework, see \cite[Chapter 3]{Bill} or \cite[Chapter 3]{WW}, for example), where $X^{SG}$ is a strong Markov diffusion -- the so-called Brownian motion on the Sierpi\'nski gasket, started from $x_0$. We remark that the terminology `Brownian motion' reflects the fact that $X^{SG}$ is apparently the most natural stochastic process on the Sierpi\'nski gasket -- apart from being a strong Markov diffusion that arises as a scaling limit of random walks on approximating lattices, it has a distribution that is invariant under the symmetries of the underlying space, and also satisfies natural scale invariance properties. Given this, as in other settings, it is natural to ask how robust a result such as \eqref{rwscaling} is to perturbations in the environment in which the process $X^n$ is based.

One simple, canonical way in which to introduce disorder into the situation is in terms of the random conductance model. Specifically, let $G=(V_G,E_G)$ be a locally finite, connected graph. Let $\omega=(\omega_e)_{e\in E_G}$ be a collection of independent and identically distributed (i.i.d.) strictly-positive random variables built on a probability space with probability measure $\mathbf{P}$; these are the so-called \emph{random conductances}. (Actually, for our model of self-similar fractals, we will later allow some local dependence.) Conditional on $\omega$, we define the \emph{variable speed random walk (VSRW)} $X^{V}=(X^{V}_t)_{t\geq 0}$ to be the continuous-time $V_G$-valued Markov chain with jump rate from $x$ to $y$ given by $\omega_{xy}$ if $\{x,y\}\in E_G$, and jump rate 0 otherwise. We obtain the associated \emph{constant speed random walk (CSRW)} $X^{C}=(X^{C}_t)_{t\geq 0}$ by setting the jump rate along edge $x$ to $y$ to be $\omega_{xy}/\nu(\{x\})$, where
 \begin{equation}\label{nudef}
 \nu\left(\{x\}\right):=\sum_{e\in E_G:\:x\in e}\omega_e;
 \end{equation}
note that the latter process has unit mean holding times at each vertex, and so $X^n$ as described in the previous paragraph is simply the CSRW when $G_n$ is equipped with unit conductances $\omega_e\equiv 1$.

An important observation is that the VSRW and CSRW experience different trapping behaviour on edges of large conductance. In particular, if we have an edge of conductance $\omega_e\gg1$ (surrounded by other edges of conductance close to 1), then both the VSRW and CSRW cross the edge order $\omega_e$ times before escaping. However, each crossing only takes the VSRW a time of $1/\omega_e$, meaning that it is only trapped for a time of order 1, whereas each crossing for the CSRW takes a time of order 1, and so the latter process is trapped for a total time of order $\omega_e$. In particular, when the weights are bounded away from 0, but not bounded above, we might expect the VSRW of the random conductance model to behave like the VSRW on the unweighted graph. For the CSRW, however, we would expect the trapping to be more significant, potentially leading to anomalous scaling if the weights are suitably inhomogeneous.

The random conductance model has been studied in a range of settings, via which the intuition of the previous paragraph has been shown to reflect the actual behaviour of the VSRW and CSRW. In the case of $\mathbb{Z}^d$ with $d\geq 2$, for example, it has been established that if the weights are bounded away from 0, then the VSRW always scales diffusively to a Brownian motion \cite{BD}. On the other hand, for the CSRW this is only true when the weights also have a finite first moment \cite{BD}. (In fact, both these results also apply when $d=1$, cf.\ remarks in \cite{Bisp,CHK}. See also \cite{ABDH} for the case when the weights are unbounded below, and \cite{ADS} for results beyond the case of i.i.d.\ conductances.) For weights whose tail no longer has a first moment, but is in the normal domain of attraction of an $\alpha$-stable random variable, namely there exists a constant $c\in(0,\infty)$ such that
\begin{equation}\label{utd}
u^{\alpha}\mathbf{P}\left(\omega_e>u\right)\rightarrow c
\end{equation}
as $u\rightarrow\infty$, one instead sees as a scaling limit for the CSRW the fractional kinetics process -- this is a Brownian motion subordinated by an $\alpha$-stable process, which is subdiffusive \cite{BCJ, CJ}. The subordination here reflects that in its first $n$ jumps, the random walk visits $Cn$ sites, and the time spent in these grows like a sum of $n$ i.i.d.\ $\alpha$-stable random variables, so is of order $n^{1/\alpha}\gg n$ (there are logarithmic corrections needed when $d=2$ \cite{CJ}). In $d=1$, the simple random walk revisits sites more often, and so although it is also true that the CSRW is subdiffusive when the weights satisfy \eqref{utd}, the nature of the process is different. Rather, the limiting process, is a Brownian motion time-changed by the Poisson random measure
\begin{equation}\label{prm}
\nu(dx)=\sum_{i}v_i\delta_{x_i}(dx),
\end{equation}
where $(v_i ,x_i)_{i\in\mathbb{N}}$ is a Poisson point process with intensity $\alpha v^{-1-\alpha}dvdx$, and $\delta_{x_i}$ is the probability measure placing all its mass at $x_i$; this random measure can be viewed as the scaling limit of the random trapping environment \cite{CJ}. After its introduction in \cite{FIN} as a scaling limit for a random walk with strongly inhomogeneous random jump rates, the Brownian motion time-changed by $\nu$ is called the Fontes-Isopi-Newman diffusion.

For fractals, the random conductance model has previously been studied in \cite{Kum2,kk}, where homogenisation was shown for certain classes of fractal graphs when the weights were bounded uniformly below and above. Here, we explain the progress of \cite{CHK}, in which a framework was developed that allowed unbounded weights, and particularly weights satisfying \eqref{utd} to be considered. For the particular case of nested fractals (the precise definition of which is recalled in the next section), one knows that diffusions on such spaces are point recurrent, and so it is natural to conjecture that the nature of the random conductance model is likely to be more closely related to the one-dimensional Euclidean picture than the higher dimensional situation. The aim of this article is to explain that this is indeed the case, with the main result being stated as Theorem \ref{ssfrcmresult}. We note that, although we restrict to nested fractals here, in \cite{CHK}, the slightly more general setting of uniformly finitely ramified fractals was considered. Moreover, we also remark that heat kernel estimates for the limiting processes are given in \cite{CroHamKum}.

The remainder of the article is organised as follows. After introducing the model in Section \ref{nfsec}, we go on to study the renormalisation and homogenisation of associated resistance metrics in Section \ref{homogsec}, and then present the main scaling result in Section \ref{concsec}.

\section{Random conductance model on nested fractal graphs}\label{nfsec}

In this section, we introduce precisely the model that will be of interest in the remainder of the article, starting with the notion of a nested fractal. For $\beta>1$ and $I=\{1,2,\cdots, N\}$, let $(\psi_i)_{i\in I}$ be a family of contraction maps on ${\mathbb{R}}^d$ such that $\psi_i(x)=\beta^{-1}U_ix +\gamma_i$ for $x\in {\mathbb{R}}^d$, where $U_i$ is a unitary map and $\gamma_i\in {\mathbb{R}}^d$. As $(\psi_i)_{i\in I}$ is a family of contraction maps, there exists a unique non-void compact set $F$ such that ${F} =\cup_{i\in I}\psi_i ({F})$. We assume the following.
\begin{description}
  \item[\textbf{Open set condition}] There is a non-empty, bounded open set $W$ such that the sets $(\psi_i (W))_{i\in I}$ are disjoint and $\cup_{i\in I} \psi_i (W)\subseteq W$.
\end{description}
The maps $(\psi_i)_{i\in I}$ have unique fixed points, and we denote the set of these by $Fix$. A point $x \in Fix$ is called an {\it essential fixed point} if there exist $i,j \in I,~i \ne j$ and $y\in Fix$ such that $\psi_{i}(x)=\psi_{j}(y)$. We write ${V}_0$ for the set of essential fixed points. Denoting $\psi_{i_{1},\dots,i_{n}}=\psi_{i_{1}}\circ\dots\circ\psi_{i_{n}}$ for each $n\ge 0$ and $i_1,\cdots,i_{n}\in I$, we call a set of the form $\psi_{i_1,\cdots, i_n}({ V}_0)$ an $n$-cell. The further assumptions we make are the following.
\begin{description}
  \item[\textbf{Connectivity}] For any $1$-cells $C$ and $C'$, there is a sequence $C=C_0,C_1,\dots,C_n=C'$ of $1$-cells such that $C_{i-1}\cap C_i\neq \emptyset$ for $i=1,\dots,n$.
  \item[\textbf{Symmetry}] For any $x,y\in\mathbb{R}^d$ with $x\neq y$, let $H_{xy}$ denote the hyperplane perpendicularly bisecting $x$ and $y$, and $U_{xy}$ denote reflection with respect to $H_{xy}$. If $x,y\in V_0$ and $x\neq y$, then $U_{xy}$ maps $n$-cells to $n$-cells, and maps any $n$-cell which contains elements on both sides of $H_{xy}$ to itself for each $n\geq 0$.
  \item[\textbf{Nesting/Finite ramification}] If $n\geq 1$ and if $(i_1,\cdots,i_{n})$ and $(j_1,\cdots,j_{n})$ are distinct elements of $I^n$, then
  \[\psi_{i_{1},\dots,i_{n}}({ F}) \bigcap\psi_{j_{1},\dots,j_{n}}
({F})=\psi_{i_{1},\dots,i_{n}}({V}_0)\bigcap\psi_{j_{1},
\dots,j_{n}}({ V}_0).\]
  \end{description}
A \emph{nested fractal} ${F}$ is a set determined by $(\psi_i)_{i\in I}$ satisfying  the above assumptions with $|{V}_0|\ge 2$. Throughout, we assume without loss of generality that $\psi_1(x)=\beta^{-1}x$ and $0$ belongs to ${V}_0$. We observe that the class of nested fractals was introduced in \cite{lind}, and is included in the class of uniformly finitely ramified fractals, first introduced in \cite{HK} (and upon which the random conductance model was studied in \cite{CHK}), and the latter collection is included in the class of post-critically finite self-similar sets \cite{kig1}.  We note that the Sierpi\'nski gasket is a nested fractal, other examples include the Vicsek set, and Lindstr{\o}m's snowflake. Some discussion about the restrictiveness of the axioms for nested fractals appears in \cite[Remark 5.25]{Barlow}.

Related to the nested fractal itself, we now introduce a sequence of nested fractal graphs $(G_n)_{n\geq 0}$. As in the case of the Sierpi\'nski gasket described in the introduction, the $G_n$ has vertex set $V_n$ given by $\cup_{i=1}^N\psi_i(V_{n-1})$, where $V_0$ is as defined above. Moreover, for each $n$, the edge set $E_n$ of $G_n$ consists of the collection of pairs of vertices that are contained in the same $n$-cell. We let $\mu_n$ be the counting measure on $V_n$ (placing mass one on each vertex).

Finally for this section, let us describe the version of the random conductance model that is of interest here. For each $n\geq 1$, let $\omega^n=(\omega^n_e)_{e\in E_n}$ be a collection of strictly-positive random variables built on a probability space with probability measure $\mathbf{P}$. We assume the following conditions on the weights.
\begin{description}
  \item[\textbf{Independence}] Weights within each $n$-cell are independent copies of $\omega^0$.
  \item[\textbf{Uniform lower bound}] There exists a deterministic constant $c>0$ such that, $\mathbf{P}$-a.s.,
  \[\omega^0_e\geq c.\]
  \item[\textbf{$\alpha$-stable tail decay}] There exist constants $\alpha\in(0,1)$ and $c\in(0,\infty)$ such that the random conductance distribution satisfies
      \begin{equation}\label{alpha}
u^{\alpha}\mathbf{P}\left(\sum_{e\in E_0}\omega^0_e>u\right)
\rightarrow c 
\end{equation}
as $u\rightarrow\infty$.
\end{description}
Given a realisation of weights satisfying these assumptions, we define the variable speed random walk $X^{n,V}$ and constant speed random walk $X^{n,C}$ on $G_n$, as per the conventions in the introduction. Specifically, both have jump chains given by the simple random walk on the graph $G_n$. The process $X^{n,V}$ has exponential holding times, with the mean of the holding times at vertex $x\in V_n$ being given by $1/\nu_n(\{x\})$, where, similarly to \eqref{nudef},
\begin{equation}\label{nudefn}
 \nu_n\left(\{x\}\right):=\sum_{e\in E_n:\:x\in e}\omega^n_e;
 \end{equation}
the process $X^{n,C}$ has unit mean exponential holding times. The so-called quenched, i.e.\ conditional on the conductances, laws of $X^{n,V}$ and $X^{n,C}$ started from a vertex $x\in V_n$ will be denoted $P^{n,V}_x$ and $P^{n,C}_x$, respectively. The corresponding averaged/annealed laws are then given by
\[\mathbb{P}^{n,V}_x:=\int P^{n,V}_x\left(\cdot\right)d\mathbf{P},\qquad\mathbb{P}^{n,C}_x:=\int P^{n,C}_x\left(\cdot\right)d\mathbf{P}.\]
The aim of this article is to describe scaling limits for both  $X^{n,V}$ and $X^{n,C}$ under their annealed laws; the main result is stated as Theorem \ref{ssfrcmresult}. Some discussion as to why we consider the annealed laws, rather than the quenched laws, is given in Remark \ref{queann}.

\section{Homogenisation of resistance}\label{homogsec}

In this section, we will briefly recall the now classical construction of a resistance metric on a nested fractal via graphical approximations. Following this, we explain what is perhaps the main result of \cite{CHK} concerning self-similar fractals, which is that the same resistance metric arises from the random conductance model defined in the previous section, i.e.\ homogenisation of the resistance occurs. Roughly speaking this can be interpreted as meaning that, apart from normalisation by a deterministic constant, the randomness of the conductances is insignificant on large scales. Intuitively, this might be expected since, whilst the tail decay at \eqref{alpha} leads to the occasional exceptionally large edge conductance, or equivalently the occasional exceptionally small edge resistance, as we rescale, neighbouring points are anyway close in terms of resistance, and so this does not lead to large scale distortions.

Before getting to resistance metrics, however, we introduce the canonical Dirichlet form and Brownian motion on a nested fractal. In Lindstr{\o}m's original work on nested fractals \cite{lind}, transition probabilities $(q_{x,y})_{x,y\in V_0}$ satisfying $q_{x,x}=0$ and $\sum_{y\in V_0}q_{x,y}=1$ for $x\in V_0$, and also $q_{x,y}=q_{y,x}>0$ for $x\neq y\in V_0$ were introduced. Importantly, it was further established that the quantities $(q_{x,y})_{x,y\in V_0}$ could be chosen to be invariant under renormalisation in the sense we now describe. Specifically, define a quadratic form by setting
\[\mathcal{E}_0(f,f)=\frac{1}{2}\sum_{x,y\in V_0}q_{x,y}\left(f(x)-f(y)\right)^2\]
for $f\in \mathcal{F}_0:=\{f:V_0\rightarrow\mathbb{R}\}$. One obtains a further quadratic form on the same space by defining
\[\tilde{\mathcal{E}}_0(f,f)=\inf\left\{\sum_{i\in I}\mathcal{E}_0\left(g\circ\psi_i,g\circ\psi_i\right):\:g:V_1\rightarrow \mathbb{R},\:g|_{V_0}=f\right\}\]
for $f\in \mathcal{F}_0$. The invariance under renormalisation of \cite[Theorem V.5]{lind} then has the equivalent statement that there exists a constant $\rho>1$ such that $\mathcal{E}_0=\rho\tilde{\mathcal{E}}_0$. Moreover, it is now known that the latter condition, together with the assumption that $q$ are the entries of a stochastic matrix, ensure the uniqueness of $(q_{x,y})_{x,y\in V_0}$ (see \cite[Theorem 6.8]{Sabot} and \cite[Corollary 3.5]{kk}). Given $(q_{x,y})_{x,y\in V_0}$ and $\rho$, for $n\geq 1$ we then let
\[\mathcal{E}_n(f,f)=\rho^n\sum_{i_{1},\dots,i_{n}\in I}\mathcal{E}_0\left(f\circ\psi_{i_{1},\dots,i_{n}},f\circ\psi_{i_{1},\dots,i_{n}}\right)\]
for  $f\in \mathcal{F}_n:=\{f:V_n\rightarrow\mathbb{R}\}$. One then obtains a canonical quadratic form on $F$ by setting
\[\mathcal{E}(f,f):=\lim_{n\rightarrow\infty}\mathcal{E}_n(f|_{V_n},f|_{V_n})\]
for any $f\in\mathcal{F}:=\{f\in C(F,\mathbb{R}):\:\lim_{n\rightarrow\infty}\mathcal{E}_n(f|_{V_n},f|_{V_n})<\infty\}$. Importantly, the resulting quadratic form $(\mathcal{E},\mathcal{F})$ turns out to be a Dirichlet form on $L^2(F,\mu)$, where $\mu$ is the unique self-similar probability measure on $F$, that is, the only probability measure satisfying
\[\mu=\frac{1}{N}\sum_{i\in I}\mu\circ \psi_i^{-1}.\]
As a consequence, standard machinery from probability theory (see \cite{FOT}, for example) yields that there exists a corresponding Markov process $X^F=(X^F_t)_{t\geq 0}$, which is now commonly called the Brownian motion on the nested fractal $F$.

We next describe the parallel construction of the resistance metric on $F$. To start with one possible definition, we observe that from the quadratic form $(\mathcal{E},\mathcal{F})$, one obtains a metric on $F$ by defining
\begin{equation}\label{rdef}
R(x,y)^{-1}:=\inf\left\{\mathcal{E}(f,f):\:f\in\mathcal{F},\:f(x)=1,\: f(y)=0\right\},\qquad x,y\in F,\: x\neq y;
\end{equation}
this is the resistance metric on $F$. In fact, the above description of $R$ yields a one-to-one relationship between a class of quadratic forms called resistance forms (of which $(\mathcal{E},\mathcal{F})$ is one) and a class of metrics called resistance metrics (see \cite[Theorems 2.3.4, 2.3.6]{kig1}, for example). An alternative definition of $R$ is via resistance metrics on the finite graphs. Specifically, suppose $R_n$ is the resistance metric on $V_n$ induced by placing conductances according to $(\rho^{-n}q_{x,y})_{x,y\in V_0}$ along edges of $n$-cells, i.e.\ setting the conductance from $\psi_{i_1,\dots,i_n}(x)$ to $\psi_{i_1,\dots,i_n}(y)$ to be  $\rho^{-n}q_{x,y}$; alternatively, $R_n$ can be defined from $(\mathcal{E}_n,\mathcal{F}_n)$ analogously to \eqref{rdef}. From the invariance under renormalisation of $\mathcal{E}_0$, one can check that
\[R_n=R_m|_{V_n}, \qquad \forall m\geq n.\]
From this it readily follows that we have $R=\lim_{n\rightarrow\infty}R_n(x,y)$ on $V_*=\cup_{n\geq 0}V_n$. In particular, $R|_{V_n}=R_n$. With some additional work to check that $(F,R)$ is the completion of $(V_*,R)$, we obtain that $V_n$ converges to $F$ with respect to Hausdorff topology on compact subsets of $(F,R)$. (See \cite{Keff} for proofs of these claims.)

It transpires that one obtains the limit described in the preceding paragraph if the deterministic conductances characterised by $(q_{x,y})_{x,y\in V_0}$ are replaced by the random conductances of the previous section. That is, suppose $R_n^\omega$ is the resistance metric on $V_n$ induced by placing conductances according to $(c\rho^{-n}\omega_e^n)_{e\in E_n}$ along edges of the graph, where $c\in(0,\infty)$ is a deterministic constant that depends on the law of the conductances; this is the metric given by \eqref{rdef} for the following quadratic form
\[\frac{1}{2c}\rho^n\sum_{i_{1},\dots,i_{n}\in I}\sum_{x,y\in V_0}\omega^n_{\psi_{i_{1},\dots,i_{n}}(x),\psi_{i_{1},\dots,i_{n}}(y)}\left(f\circ\psi_{i_{1},\dots,i_{n}}(x)-f\circ\psi_{i_{1},\dots,i_{n}}(y)\right)^2,\]
which is defined for $f\in \mathcal{F}_n$. From \cite[Theorem 6.11]{CHK}, we then have that, in $\mathbf{P}$-probability,
\begin{equation}\label{rconv}
\left(R^\omega_n(x,y)\right)_{x,y\in V_0}\rightarrow \left(R(x,y)\right)_{x,y\in V_0},
\end{equation}
where we note that the constant $c$ is determined by this result. The proof in \cite{CHK}, which can heuristically be understood as establishing contractivity of a renormalisation map, resembles that of the corresponding results in \cite{Kum2,kk}. However, the lack of a uniform upper bound on the conductances leads to significant technical challenges, particularly in checking that certain quantities are integrable, as is required for the argument to work. From \eqref{rconv} and the trivial bound that $R_n^\omega\leq CR_n$, (which follows from the fact that the conductances are bounded away from 0,) we readily obtain the following proposition.

{\proposition[{\cite[Lemma 6.14]{CHK}}] \label{ssfresconv} In $\mathbf{P}$-probability,
\[\sup_{x,y\in V_n}\left|R_n^\omega(x,y)-R(x,y)\right|\rightarrow 0.\]}

Since $(V_n,R_n^\omega)$ can not in general be isometrically embedded into $(F,R)$, then the usual Hausdorff topology on $(F,R)$ is not the right topology with which to discuss convergence. However, one can instead conclude from the previous result (and some small additional technical work again depending on the bound $R_n^\omega\leq CR_n$) that $(V_n,R_n^\omega)$ converges to $(F,R)$ with respect to the Gromov-Hausdorff topology, that is, all the spaces in question can be isometrically embedded into a common metric space so that the $V_n$ converges to $F$ with respect to the usual Hausdorff metric on this space (see \cite[Chapter 7]{BBI} for background on the Gromov-Hausdorff topology).

\section{Random walk scaling limits}\label{concsec}

Proposition \ref{ssfresconv} is the main ingredient to proving scaling limits for the variable speed random walk $X^{n,V}$ and  the constant speed random walk $X^{n,C}$. Indeed, the only additional input required is the convergence under scaling of the counting measure $\mu_n$ and the measure $\nu_n$ defined in terms of conductances at \eqref{nudefn}, which is straightforward to prove. The machinery that allows us to proceed with this program is the main result of \cite{C} (which gives a more general version of the result of \cite{CHK}).

To introduce the abstract result we appeal to precisely, let us fix the framework. In particular, we write $\mathbb{F}^*_c$ for the collection of quintuples of the form $(K,R_K,\mu_K,\rho_K,\phi_K)$, where: $K$ is a non-empty set; $R_K$ is a resistance metric on $K$ such that $(K,R_K)$ is compact; $\mu_K$ is a locally finite Borel regular measure of full support on $(K,R_K)$; $\rho_K$ is a marked point in $K$, and $\phi_K$ is a continuous map from $K$ to some fixed metric space $(M,d_M)$. From the point of view of metric geometry, there is a natural notion of convergence of such spaces which gives rise to the marked spatial Gromov-Hausdorff-Prohorov topology. Specifically, convergence of some sequence in $\mathbb{F}^*_c$ means that all the spaces can be isometrically embedded into a common metric space $(\mathcal{M},d_\mathcal{M})$ in such a way that: the embedded sets converge with respect to the Hausdorff distance, the embedded measures converge weakly, the embedded marked points converge, and the image of the continuous map is close in $M$ for points that are close in $\mathcal{M}$. We note that such Gromov-Hausdorff-type topologies have proved useful for studying various kinds of random metric spaces; see \cite{BBI} for an introduction to the classical theory. More specifically, the marked spatial Gromov-Hausdorff-Prohorov topology was introduced in \cite{BCK}, building on the notions of the Gromov-Hausdorff-Prohorov/Gromov-Hausdorff-vague topologies of \cite{ADH, ALWtop, Evans, Miermont} and the topology for spatial trees of \cite{DL} (cf.~the spectral Gromov-Hausdorff topology of \cite{CHK}).

Importantly, that the elements $(K,R_K,\mu_K,\rho_K,\phi_K)$ of $\mathbb{F}^*_c$ incorporate a resistance metric means that there is a naturally associated stochastic process. For, it is a result of Kigami that the corresponding resistance form, characterised via \eqref{rdef}, is a regular Dirichlet form on $L^2(K,\mu_K)$, and so naturally associated with a Markov process (see \cite{Kig}, Chapter 9, for example). The following result establishes that, if the convergence described in the previous paragraph occurs, then we also obtain convergence of stochastic processes.

\begin{theorem}[{\cite[Theorem 7.2]{C}}]\label{mainres} Suppose that $(K_n,R_{K_n},\mu_{K_n},\rho_{K_n},\phi_{K_n})_{n\geq 1}$ is a sequence in $\mathbb{F}^*_c$ satisfying
\begin{equation}\label{ghp}
\left(K_n,R_{K_n},\mu_{K_n},\rho_{K_n},\phi_{K_n}\right)\rightarrow \left(K,R_K,\mu_K,\rho_K,\phi_K\right)
\end{equation}
in the marked spatial Gromov-Hausdorff-Prohorov topology for some element $(K,R_K,\mu_K,\rho_K,\phi_K)\in \mathbb{F}^*_c$. It then holds that
\[P^n_{\rho_{K_n}}\left(\left(\phi_{K_n}(X^n_t)\right)_{t\geq 0}\in\cdot\right)\rightarrow P_{\rho_K}\left(\left(\phi_K(X_t)\right)_{t\geq 0}\in\cdot\right)\]
weakly as probability measures on $D(\mathbb{R}_+,M)$, where  $((X^n_t)_{t\geq 0},(P^n_x)_{x\in K_n})$ is the Markov process corresponding to $(K_n,R_{K_n},\mu_{K_n},\rho_{K_n})$, and $((X_t)_{t\geq 0},(P_x)_{x\in K})$ is the Markov process corresponding to $(K,R_K,\mu_K,\rho_K)$.
\end{theorem}

{\remark The key to the proof of the above result in \cite{C} is the observation that for a process associated with a resistance metric, it is possible to explicitly express the associated resolvent kernel in terms of the resistance metric. (This was also the basis of the corresponding argument for trees from \cite{ALWtree}.) Specifically, if $((X_t)_{t\geq 0},(P_x)_{x\in K})$ is the Markov process associated with $(K,R_K,\mu_K,\rho_K,\phi_K)\in\mathbb{F}^*_c$, define the resolvent of $X$ killed on hitting $x$ by
\[G_xf(y)=E_y\int_0^{\sigma_x}f(X_s)ds,\]
where $E_y$ is the expectation under $P_y$, and $\sigma_x:=\inf\{t\geq 0:\:X_t=x\}$ is the hitting time of $x$ by $X$. (NB.\ Processes associated with resistance forms hit points; the above expression is well-defined and finite.) One can then write
\[G_xf(y)=\int_Kg_x(y,z)f(z)\mu_K(dz),\]
where the resolvent kernel is given by
\[g_x(y,z)=\frac{R_K(x,y)+R_K(x,z)-R_K(y,z)}{2}.\]
(See \cite[Theorem 4.3]{Kig}.) Appealing to this formula, the metric measure convergence at (\ref{ghp}) enables one to check the convergence of resolvents in a certain sense. One can then use more standard machinery from probability theory to establish semigroup convergence, and moreover convergence of finite dimensional distributions. To complete the proof, one is also required to check tightness of the processes (see \cite[Chapter 16]{Bill}), but again this can be deduced from the above resolvent density formula (or, more precisely, a slight generalisation thereof). See \cite{C} for details.}

{\remark Whilst Theorem \ref{mainres} has an appealingly concise statement, checking the assumption at \eqref{ghp} is by no means trivial. Indeed, beyond the case of graph trees (or graphs that are close to trees), where the resistance metric corresponds to (or is close to, respectively) a shortest path metric, or certain finitely ramified self-similar fractals, where the resistance metric can be studied by using the particular structure of the space, understanding detailed properties of the resistance metric remains a challenge. To give just one example of an open problem from the world of self-similar fractals, it is still not known how to compute the value of the resistance exponent for graphs based on the two-dimensional Sierpi\'nski carpet, see \cite{BBres} for some work in this direction, and the discussion in \cite[Example 4]{BCoulK} concerning the graphical Sierpi\'nski carpet in particular.}
\bigskip

We will apply Theorem \ref{mainres} with $K_n=V_n$, $R_{K_n}=R_n^\omega$, $\mu_{K_n}=\mu_n$ or $\mu_{K_n}=\nu_n$, $\rho_{K_n}=0$, and $\phi_{K_n}:=I_n$, where $I_n$ is the identity map from $K_n$ into $\mathbb{R}^d$. The following lemma gives us the scaling limits of the measures. To state the result, we introduce a Poisson random measure on $F$ by setting
\[\nu(dx)=\sum_{i}v_i\delta_{x_i}(dx),\]
where $(v_i ,x_i)_{i\in\mathbb{N}}$ is a Poisson point process with intensity $\alpha v^{-1-\alpha}dv\mu(dx)$, and $\delta_{x_i}$ is the probability measure placing all its mass at $x_i$. (This is the analogue of the measure defined at \eqref{prm} in the present setting.) Note that the exponent $\alpha$ is given by the tail of the conductance distribution \eqref{alpha}.

{\lemma\label{mconv} It holds that $N^{-n}\mu_n\rightarrow \mu$, and also there exists a deterministic constant $c_0\in(0,\infty)$ such that $c_0^{-1}N^{-n/\alpha}\nu_n\rightarrow \nu$ in distribution, in both cases with respect to the weak topology for finite measures on $\mathbb{R}^d$.}
\bigskip

Combining Proposition \ref{ssfresconv} and Lemma \ref{mconv}, we readily obtain that
\begin{equation}\label{wlln}
\left(V_n,R^\omega_n,N^{-n}\mu_n,0,I_n\right)\rightarrow \left(F,R,\mu,0,I\right),
\end{equation}
in $\mathbf{P}$-probability, and
\[\left(V_n,R^\omega_n,c_0^{-1}N^{-n/\alpha}\nu_n,0,I_n\right)\rightarrow \left(F,R,\nu,0,I\right),\]
in distribution under $\mathbf{P}$ with respect to the marked spatial Gromov-Hausdorff-Prohorov topology, where $I$ is the identity map from $F$ into $\mathbb{R}^d$. Since $X^{n,V}$ is the process associated with $(V_n,c^{-1}\rho^nR^\omega_n,\mu_n,0,I_n)$, and $X^{n,C}$ is the process naturally associated with $(V_n,c^{-1}\rho^nR^\omega_n,\nu_n,0,I_n)$, we are consequently able to apply Theorem \ref{mainres} to deduce a scaling limit for these processes. (By considering the generators of the relevant Markov processes, it is readily checked how the resistance and mass scaling factors can be interpreted in terms of time scaling.) As for the limiting processes, we note that the Brownian motion $X^F$ is the process associated with $(F,R,\mu,0)$ -- we write the law of this process started from 0 as $P_0$. Moreover, the process associated with $(F,R,\nu,0)$ is the time-change of $X^F$ according to $\nu$, that is, defining an additive functional
\[A_t:=\int_0^tL_t(x)\nu(dx),\]
where $(L_t(x))_{x\in F,\:t>0}$ are the jointly continuous local times of $X^F$ (with respect to $\mu$), and its right-continuous inverse $\tau(t):=\inf\{s>0:\:A_s>t\}$, we set
\[X^{F,\nu}_t:=X^F_{\tau(t)};\]
following the definition of the corresponding one-dimensional process in \cite{FIN}, we call this the FIN diffusion on $F$. The averaged/annealed law of the FIN diffusion on $F$, started from 0, will be denoted
\[\mathbb{P}^{{\rm FIN}}_0:=\int P_0\left(X^{F,\nu}\in\cdot\right)d\mathbf{P},\]
i.e.\ one chooses $\nu$ according to $\mathbf{P}$, and then the law of $X^{F,\nu}$ is determined by the law of $X^F$ under $P_0$.

{\theorem\label{ssfrcmresult} There exists a deterministic constant $c_1\in(0,\infty)$ such that
\[\mathbb{P}^{n,V}_{0}\left(\left(X^{n,V}_{c_1t(\rho N)^n}\right)_{t\geq 0}\in\cdot\right)\rightarrow {P}_0\left(\left(X^F_{t}\right)_{t\geq 0}\in\cdot\right)\]
weakly as probability measures on $D(\mathbb{R}_+, \mathbb{R}^d)$. Moreover, there exists a deterministic constant $c_2\in(0,\infty)$ such that
\[\mathbb{P}^{n,C}_{0}\left(\left(X^{n,C}_{c_2 t (\rho N^{1/\alpha})^n}\right)_{t\geq 0}\in\cdot\right)
\rightarrow\mathbb{P}^{{\rm FIN}}_{\rho}\left(\left(X^{F,\nu}_{t}\right)_{t\geq 0}\in\cdot\right)\]
weakly as probability measures on $D(\mathbb{R}_+,\mathbb{R}^d)$.}

{\remark To state the result for the Sierpi\'nski gasket explicitly, note that in this case we have $N=3$ and $\rho=5/3$, so that
\[\mathbb{P}^{n,V}_{0}\left(\left(X^{n,V}_{c_1t5^n}\right)_{t\geq 0}\in\cdot\right)\rightarrow {P}_0\left(\left(X^F_{t}\right)_{t\geq 0}\in\cdot\right),\]
and we also have
\[\mathbb{P}^{n,C}_{0}\left(\left(X^{n,C}_{c_2 t5^n(3^{\frac{1}{\alpha}-1})^n}\right)_{t\geq 0}\in\cdot\right)
\rightarrow\mathbb{P}^{{\rm FIN}}_{0}\left(\left(X^{F,\nu}_{t}\right)_{t\geq 0}\in\cdot\right).\]
In particular, the scaling regime for the variable speed random walk matches that of the simple random walks on the unweighted graphs, as stated at \eqref{rwscaling}; and since $\alpha<1$, the constant speed random walk (or limiting diffusion) moves through the relevant graph more slowly than the unweighted simple random walk (or Brownian motion, respectively). Together with known results for simple random walks on nested fractal graphs, Theorem \ref{ssfrcmresult} implies that these qualitative comments apply to nested fractal graphs in general.}

{\remark When $\mathbf{E}\omega_e^0<\infty$ for each $e\in E_0$, one obtains in place of the second claim of Lemma \ref{mconv} that there exists a constant $c_0$ such that $c_0^{-1}N^{-n}\nu_n\rightarrow\mu$. Consequently, if \eqref{alpha} is replaced by the assumption of finite first moments, then one can check the annealed limit of $X^{n,C}$ is Brownian motion, rather than the FIN diffusion that appears in the second statement of Theorem \ref{ssfrcmresult}.}

{\remark\label{queann} A stronger notion of convergence than convergence with respect to the annealed law is convergence with respect to the quenched law for $\mathbf{P}$-a.e.\ realisation of the conductances. Typically, one might hope to be able to prove such a quenched convergence statement in the case where the conductances homogenise, as has been established when the underlying graph is a Euclidean lattice (see \cite{ABDH, ADS, BD}, for example). In particular, it would be natural to conjecture that for the example described in this article, the quenched law of the VSRW $X^{n,V}$ converges as $n\rightarrow\infty$ for typical realisations of the environment. To do this, it would be sufficient to replace the weak (i.e.\ in probability) statement of \eqref{wlln} with a strong (i.e.\ $\mathbf{P}$-a.s.) one. However, the techniques of \cite{CHK} are not sufficient to yield such a result. As for the CSRW $X^{n,C}$, the typical fluctuations of the conductance environment as $n$ varies will be too large to permit a quenched limit statement (cf.\ the law of the iterated logarithm for simple random walk on $\mathbb{Z}$, which implies that individual sample paths can not be rescaled to a realisation of Brownian motion on $\mathbb{R}$, even though the discrete paths have the latter process as a distributional limit).}

\begin{acknowledgement}
The author is grateful to the organisers of the conference Fractal Geometry and Stochastics 6 for arranging a wonderful meeting, where he was given the chance to present the work of \cite{C,CHK}, and for inviting him to produce this article. His attendance at the latter event was partially supported by the JSPS Grant-in-Aid for Research Activity Start-up, 18H05832. The author is grateful to a referee for their careful reading of an earlier version of the article, which led to a number of improvements being made.
\end{acknowledgement}

\providecommand{\bysame}{\leavevmode\hbox to3em{\hrulefill}\thinspace}
\providecommand{\MR}{\relax\ifhmode\unskip\space\fi MR }
\providecommand{\MRhref}[2]{%
  \href{http://www.ams.org/mathscinet-getitem?mr=#1}{#2}
}
\providecommand{\href}[2]{#2}

\end{document}